\documentclass[12pt,a4paper]{amsart}
\usepackage{geometry}
\usepackage{amsmath}
\usepackage{amssymb}
\usepackage{amsthm}
\usepackage{graphicx}
\usepackage{color}
\usepackage{subfig}
\newtheorem{theorem}{Theorem}
\newtheorem{lemma}[theorem]{Lemma}
\newtheorem{corollary}[theorem]{Corollary}
\newtheorem{proposition}{Proposition}

\theoremstyle{definition}
\newtheorem{remark}{Remark}
\newtheorem{example}{Example}

\def\blem{\begin{lemma}}
\def\elem{\end{lemma}}
\def\bthm{\begin{theorem}}
\def\ethm{\end{theorem}}
\def\bmat{\begin{pmatrix}}
\def\emat{\end{pmatrix}}
\def\beq{\begin{equation}}
\def\eeq{\end{equation}}
\def\erf{\eqref}

\def\R{\mathbb{R}}
\def\N{\mathbb{N}}
\def\ol{\overline}
\def\fr{\frac}
\def\mid{\,:\,} 
\def\gl{\lambda}
\def\gs{\sigma}
\def\IN{\text{ in }}
\def\INT{\operatorname{int}}

\def\supp{\operatorname{supp}}
\def\dH{\operatorname{d_H}}
\def\AND{\text{ and }}
\def\FOR{\text{ for }}

\def\ON{\text{ on }}
\def\bproof{\begin{proof}}
\def\eproof{\end{proof}}
\def\stm{\setminus}

\def\pl{\partial}
\def\gd{\delta}
\def\tim{\times}
\def\bald{\begin{aligned}}
\def\eald{\end{aligned}}
\def\FORALL{\text{ for all }}
\def\cP{\mathcal{P}}
\def\bS{\mathbb{S}}

\def\blem{\begin{lemma}}
\def\elem{\end{lemma}}
\def\bthm{\begin{theorem}}
\def\ethm{\end{theorem}}
\def\bmat{\begin{pmatrix}}
\def\emat{\end{pmatrix}}
\def\beq{\begin{equation}}
\def\eeq{\end{equation}}
\def\erf{\eqref}

\def\R{\mathbb{R}}
\def\N{\mathbb{N}}
\def\ol{\overline}
\def\fr{\frac}
\def\mid{\,:\,} 
\def\gl{\lambda}
\def\ga{\alpha}
\def\IN{\text{ in }}
\def\INT{\operatorname{int}}

\def\AND{\text{ and }}
\def\FOR{\text{ for }}

\def\ON{\text{ on }}
\def\bproof{\begin{proof}}
\def\eproof{\end{proof}}
\def\stm{\setminus}

\def\pl{\partial}
\def\gd{\delta}
\def\tim{\times}
\def\bald{\begin{aligned}}
\def\eald{\end{aligned}}
\def\FORALL{\text{ for all }}
\def\IF{\text{ if }}
\def\cP{\mathcal{P}}
\def\stm{\setminus}
\def\gO{\Omega}
\def\dist{\operatorname{dist}}

\def\du#1{\langle#1\rangle}
\def\LSC{\operatorname{LSC}}
\def\C{\operatorname{C}}
\def\ep{\varepsilon}
\def\bcases{\begin{cases}}
\def\ecases{\end{cases}}

\def\cC{\mathcal{C}}
\def\cB{\mathcal{B}}
\def\bye{\end{document}}

\def\diag#1{\operatorname{diag}(#1)}

\def\gth{\theta}
\def\gk{\kappa}

\definecolor{labelkey}{rgb}{0.2,0.2,0.2}
\definecolor{gray}{gray}{0.5}

\title[Positivity sets and the strong maximum principle]{Positivity sets of 
supersolutions of degenerate elliptic equations and the strong maximum principle}

\author[I. Birindelli, G. Galise]{Isabeau Birindelli, Giulio Galise}
\address{Dipartimento di Matematica G. Castelnuovo, Sapienza Universit\`a
di Roma, P.le Aldo Moro 2, I-00185 Roma, Italy}
\email{isabeau@mat.uniroma1.it}
\email{galise@mat.uniroma1.it}

\author[H. Ishii]{Hitoshi Ishii}
\address{Institute for Mathematics and Computer Science\\ Tsuda  University\\
 2-1-1 Tsuda, Kodaira\\ Tokyo 187-8577\\Japan}
\email{hitoshi.ishii@waseda.jp}

\thanks{
The work of HI was partially supported by the JSPS grants: KAKENHI
\#16H03948, \#18H00833. HI thanks the Department of Mathematics at the Sapienza University of Rome for financial support and its hospitality while his visit there in May 6--June 5  2019.}

\begin{document}

\begin{abstract}
We investigate positivity sets of nonnegative supersolutions of the 
fully nonlinear elliptic equations $F(x,u,Du,D^2u)=0$ in $\gO$, where $\gO$ is an open subset of $\R^N$, and the validity of the strong maximum principle for 
$F(x,u,Du,D^2u)=f$ in $\gO$, with $f\in \C(\gO)$  being nonpositive.   We obtain geometric characterizations of positivity sets $\{x\in\gO\mid u(x)>0\}$
of nonnegative supersolutions $u$
and establish the strong maximum principle under some geometric 
assumption on the set $\{x\in\gO\mid f(x)=0\}$.   
\end{abstract}

\maketitle

\section{Introduction}

\noindent 
We consider fully nonlinear degenerate elliptic equations   
that do not satisfy the strong maximum principle.
More precisely we deal with second order equations of the type 
\begin{equation}\label{intro1}
F(x,u,Du,D^2u)=0\qquad\text{in \;$\gO$,}
\end{equation}
where $\gO$ is an open subset of $\R^N$,  
such that there exists  $u\in \LSC(\gO)$  a viscosity supersolution of \eqref{intro1}, verifying  the conditions
\begin{equation}\label{intro1-1}
\min_{\gO}u=0\quad\;\text{and}\;\quad U=\left\{x\in\gO\,:\,u(x)>0\right\}\neq\emptyset.
\end{equation}
	
One of the main goals of  this paper is to provide  geometric characterizations of the positivity set $U$, defined in \erf{intro1-1}, in a quite general framework, then to analyze  more specifically the following cases:
\begin{equation}\label{intro2}
	{\mathcal P}^-_k(D^2u):=\sum_{i=1}^k\lambda_i(D^2u)=0\;\quad\text{in\; $\gO$}
	\end{equation}
and
	\begin{equation}\label{intro3}
	\lambda_k(D^2u)=0	\;\quad\text{in \;$\gO$},
	\end{equation}
where $k<N$ is a positive integer and $\lambda_1(D^2u),\ldots,\lambda_N(D^2u)$ are 
the eigenvalues of the Hessian of $u$ arranged in nondecreasing order, 
$\lambda_i\leq\lambda_{i+1}$. 
	We refer to \cite{BR, CDLV, HL1,HL1.0, HL2,  V, BGI, BGI2, BGI3,  BGL} for some work related to these 
operators in the framework of viscosity solutions. 
	
In general, if $F$ is uniformly elliptic, by the strong maximum principle, 
the set $u^{-1}(0)$ of the minimum points  is $\Omega$. Instead, when the operator is degenerate elliptic, the question is relevant and goes back to the seminal papers of Bony \cite{Bo}, Hill \cite{Hill} and  Redheffer \cite{R}, where this characterization is called the \lq\lq sharp\rq\rq\ maximum principle. 

Here, differently from all previous works, we characterize the 
 set $u^{-1}(0)$ 
 through a \lq\lq viscosity\rq\rq\ condition.  That condition has an explicit geometric  meaning when the operator $F$ is either ${\mathcal P}^-_k(D^2u)$ or $\lambda_k(D^2u)$. 
To be a little more specific we prove that, for $ {\mathcal P}^-_k(D^2u)\leq 0$, the relative boundary $\Omega\cap\partial U$ satisfies 
the inequality
\begin{equation}\label{eq1intro}
\kappa_{N-k}+\cdots+\kappa_{N-1}\geq 0, 
\end{equation}
where $\kappa_i$ denote the \lq\lq principal curvatures\rq\rq\ of $\Omega\cap\partial U$. Conversely, any subset of $\Omega$  whose boundary satisfies the above inequality is the positive set of a nontrivial supersolution of ${\mathcal P}^-_k(D^2u)\leq 0$ in $\Omega$, see Theorems \ref{FU-nec}-\ref{FU-suf} and Remark \ref{rmk2}.\\
For supersolutions of $\lambda_k(D^2u)\leq0$,  the relative boundary $\Omega\cap\partial U$ is characterized by the inequality 
\begin{equation}\label{eq2intro}
\kappa_{N-k}\geq0.
\end{equation}
In the case $k=N-1$, i.e. $\kappa_{1}\geq0$, it is nothing else but the convexity of 
the connected components of the positive set $U$, see Theorem \ref{convexity_gl_N-1}.
Let us emphasize that, since $\partial U$ is in general non smooth, the  inequalities 
\eqref{eq1intro}-\eqref{eq2intro} have to be understood in a viscosity sense which is 
exactly the purpose of the conditions \erf{FU2}-\erf{FU3} that will be introduced in 
Sections \ref{OpP_k}-\ref{Opl_k}. The conditions \erf{FU2}-\erf{FU3}  are given through a 
\lq\lq local smooth test function\rq\rq\ and this is the reason why, in analogy with viscosity solutions, we 
call them viscosity conditions.
But, as in all good generalizations, they reduce to the classical sense when $\partial U$ is a smooth  hypersurface, see Theorems \ref{thcurv} and \ref{thmcurv}.
Let us mention that 
(G$_{\gl_{N-k},\gO,U}$) is proved to be equivalent to another geometrical condition denoted \erf{CkOmU} that somehow is related to the dimension of convex subsets of $U$, see Theorem \ref{thequiv}.

In the linear case, i.e. for supersolutions of
$$
{\rm tr}(A(x)D^2u)+b(x)\cdot Du +c(x)u\leq0$$
with $A(x)\geq0$, the \lq\lq diffusion\rq\rq\ of the minimum points follows the directions of subunit vector fields $Z$ for $A$ ($Z$ is a subunit vector for $A$ if $A-Z\otimes Z\geq0$, where $(Z\otimes Z)_{ij}=Z_iZ_j$). More precisely, given $x_0\in\Omega$ such that $\min_\Omega u=u(x_0)$, then the set of minimum points contains all the points that can be reached from $x_0$ following a finite number of trajectories of subunit vector fields backward and forward in time.
This result, proved in \cite{Bo, Hill, R} for  $A(x)=\frac12\sigma(x)\sigma^T(x)$ for which the  subunit vector fields are the columns of $\sigma$, has been extended by Bardi-Da Lio in \cite{BDL1} to fully nonlinear convex operators of Bellman type, then in \cite{BDL2} for concave operators within the theory of differential games.
More recently, in \cite{BG} Bardi and Goffi have obtained new results concerning strong maximum principle and propagation of minimum points by assuming the existence of generalized subunit vectors $Z$ for the operator $F$, i.e.  $Z=Z(x)\in\R^N$ such that 
$$
\sup_{\gamma>0} F(x,0,p,\gamma p\otimes p-I)>0\quad\forall p\in\R^N,\;Z\cdot p\neq0
$$ 
(note that we adopt a reversed definition of \lq\lq ellipticity\rq\rq\ with respect to \cite{BG}). 
We stress that the operator ${\mathcal P}^-_k$ and $\lambda_k$ do not admit 
non-zero subunit vectors, since a straightforward computation shows that for any $p\in\R^N$
$$\sup_{\gamma>0} {\mathcal P}^-_k(\gamma p\otimes p-I)=-k\quad\text{and}\quad\sup_{\gamma>0} \lambda_k(\gamma p\otimes p-I)=-1.$$
Hence the previous results do not apply to our cases. 
Moreover the description of the diffusion of the minimal points seems to us to be very different from the one given here.

\smallskip
In their acclaimed works \cite{HL1,HL2,IHHL} , Harvey and Lawson  give boundary convexity conditions on $\Omega$ and $F$ 
in order to prove existence and uniqueness of 
solutions for the equation  $F(x,u,D u, D^2u)=0$ in $\Omega$ with Dirichlet boundary condition.
When $F={\mathcal P}^-_k(D^2u)$ or $F=\lambda_k(D^2u)$ the condition is similar to 
that of \eqref{eq1intro} and \eqref{eq2intro}.
Even though, it is well known that, thanks to Perron's method, existence of solutions 
is proved via the construction of sub and supersolutions, 
the scope of the supersolutions given in Theorems \ref{FU-suf0} and 
\ref{FU-suf} is of a different nature. Indeed, given a subset $U$ of $\Omega$ satisfying 
\eqref{eq1intro},  we  construct a supersolution which is positive in  $U$ and identically 
zero outside. Note that the 
standard construction, see e.g. \cite{primer}, of positive supersolutions that are zero on 
the boundary of $U$ would lead to a function that is negative outside of $U$ and that, extended to zero outside of $U$, would not be a supersolution anymore.


\medskip

The other objective of this paper is to establish some conditions on the zero level set of  nonpositive functions $f\in \C(\gO)$  that ensure that  the strong maximum principle holds
for 
\beq \label{intro4}
F(x,u,Du,D^2u)=f \quad \IN \gO.
\eeq 
The strong maximum principle states that 
\[\tag{SMP}\label{SMP}\left\{\,\begin{minipage}{0.85\textwidth}
if \ $u\in\LSC(\gO)$ \ is a supersolution of \erf{intro4} and \ $u\geq 0$ \ in \ $\gO$, then either \ $u>0$ \ in \ $\gO$\ or\ $u\equiv 0$.  \end{minipage}
\right.\] 
As already mentioned, when $f\equiv 0$, the strong maximum principle does not hold in general in our framework. However, there is a situation 
where \erf{SMP} does hold due to the geometry of the set 
\[
\Gamma=\{x\in\gO\mid f(x)=0\},
\] 
which forces, in a sense, the negativity of $f$ to be large enough.  
We establish the strong maximum principle for \erf{intro4} under 
some hypotheses on the geometry of $\Gamma$.  
Our strong maximum principle is slightly different from the classical \erf{SMP} 
and its conclusion is the positivity of nonnegative supersolutions 
of \erf{intro4} without alternative claim. 

In a certain sense this result is complementary to the one of Harvey and Lawson in \cite{HL3} where they characterize the operators for which the strong maximum principle holds. Instead we consider operators for which it does not unless the forcing term is negative enough.

Anecdotally, let us mention that even if the paper is wholly concerned with the understanding of the strong maximum principle, we never use or prove a Hopf lemma.


\medskip

The paper, after this introduction, is organized as follows.  
Section \ref{posset} is devoted to establishing theorems characterizing 
the geometry of the positivity set $U$, given in \erf{intro1-1}, 
for supersolutions of \erf{intro1}.  It introduces some basic assumptions on 
the function (or operator) $F$ as well as a geometric condition on $U$, relative to 
$\gO$, called (G$_{F,\gO,U}$). In Section \ref{OpP_k}, the condition (G$_{F,\gO,U}$), with $F=\cP_k^-$,  is examined in relation with the ``truncated mean curvature'', 
$\gk_{N-k}+\cdots+\gk_{N-1}$, where the $\gk_i$ are the principal curvatures 
of $\pl U$ at least when $U$ has a smooth boundary.  
Section \ref{Opl_k} is devoted to investigating the condition (G$_{F,\gO,U}$), with $F=\gl_k$, where a new geometric condition (C$_{k,\gO,U}$) is introduced 
and proved to equivalent to  (G$_{\gl_{N-k},\gO,U}$). Also, the convexity in the condition (C$_{1,\R^N,U}$) is studied and some related examples are presented. 
In Section \ref{smp}, the strong maximum principle is established for \erf{intro4}, together with two counterexamples. In the case of operators $\gl_k$ and $\cP_k^-$, 
a couple of a little more clear conditions for the validity of the maximum principle are given.

\subsection*{Notation:} 
$B_r^m(x)$ (resp., $\ol B_r^m(x)$) denotes the open (resp., closed) ball in $\R^m$ of radius $r$ and centered at $x$. When $x$ is the origin, we write 
$B_r^m$ and $\ol B_r^m$ instead for $B_r^m(x)$ and $\ol B_r^m(x)$, respectively. 
We write $0_m$ for the zero vector in $\R^m$. 
When $m=N$, we omit writing superscript ``$N$'', so that $B_r^N(x)=B_r(x)$
and $0_N=0$ for instance.  
For $x,y\in\R^m$, $\du{x,y}$ denotes the Euclidean inner product of $x,y$. 
For $1\leq i\leq N$, $e_i$ denotes the unit vector in $\R^N$ having unity as its $i$th entry.  The collection $\{e_1,\ldots,e_N\}$ is called the standard basis of $\R^N$. 
We write $\diag{a_1,\ldots,a_m}$ for  
the $m\tim m$ diagonal matrix with the diagonal entries $a_1,\ldots,a_m$. 		
		
\section{The positivity set}\label{posset}

Let  $\mathbb S^N$ be the linear space of $N\times N$ real symmetric matrices equipped with its usual order. We introduce the hypotheses on $F:\gO\times\R\times\R^N\times\mathbb S^N\mapsto\R$:  
\begin{enumerate}
\item[(F1)] $F\in \C(\gO\tim\R\tim\R^N\tim\bS^N)$. 
\item[(F2)] $F$ is degenerate elliptic, i.e. $F(x,u,p,X)\leq F(x,u,p,Y)$ if $X\leq Y$. 
\end{enumerate}

\begin{enumerate}
\item[(F3)] $F(x,0,0,0)\leq 0$ for all $x\in\gO$.
\item[(F4)] $F$ is positively homogeneous in the following sense: 
there is a function $h:(0,\infty)\to (0,\,\infty)$ such that
\[
F(x,u,tp,tX)=h(t)F(x,u,p,X) \ \ \FORALL t>0.
\]
\item[(F5)] $F=F(x,u,p,X)$ is \lq\lq nonincreasing\rq\rq\ in $u\in\R$. More precisely, 
$F(x,u,p,X)\leq F(x,0,p,X)$ if $u\geq 0$. 
\end{enumerate}
\medskip


\medskip
		Let $U$ be an open subset of $\gO$.
We introduce the geometric condition that 
\beq\tag{G$_{F,\gO,U}$}\label{FOmU}
\left\{\,\begin{minipage}{0.85\textwidth}
if \ $\phi\in C^2(\gO)$,\ $\hat x\in \gO\stm U$
\ and \ $\phi\leq \phi(\hat x)=0$ in $\gO\stm U$, 
then 
\[
F(\hat x, \phi(\hat x),D\phi(\hat x),D^2\phi(\hat x))\leq 0.
\]
\end{minipage}
\right. \eeq
When \erf{FOmU} holds, we say as well that 
$U$ (resp., $(\gO,U)$) satisfies (G$_{F,\gO}$) (resp., (G$_F$)). 
When $\gO=\R^N$, we refer to the conditions (G$_{F,\gO,U}$) and (G$_{F,\gO}$) above as {(G$_{F,U}$) and (G$_F$), respectively.  
\medskip

In the following two theorems we show  that the positivity set of any supersolution of \eqref{intro1} has the property \erf{FOmU} and conversely, under a slightly different assumptions on $F$ that given an open set $U$ with \erf{FOmU} property, there exists a supersolution of \eqref{intro1}  having $U$ as positivity set.

\bthm \label{FU-nec} Assume \emph{(F1), (F2),} and \emph{(F4).} 
Let $u\in\LSC(\gO)$ be a supersolution of 
\[\tag{\ref{intro1}}
F(x,u,Du,D^2u)=0 \ \ \IN \gO.
\]
Assume that $\min_{\gO}u=0$. Set 
$U=\{x\in\gO\mid u(x)>0\}$.  Then 
$U$ satisfies \emph{(G$_{F,\gO}$)}. 
\ethm

It should be noticed that, due to the lower semicontinuity of $u$, the set $U$ in the theorem above is an open set.

\bproof Let $\phi\in C^2(\gO)$ and $\hat x\in \gO\stm U$. Assume that
$\phi\leq 0=\phi(\hat x)$ in $\gO\stm U$. We may moreover assume, if needed, by replacing $\phi$ by a new one without changing the value 
$(\phi,D\phi,D^2\phi)$ at $\hat x$ that $\phi(x)<0$ for all $x\in \gO\stm (U \cup\{\hat x\})$. 
Fix $r>0$ so that $\ol B_r(\hat x)\subset\gO$. 
For any $n\in\N$ we consider the function 
$nu-\phi$ on $\ol B_r(\hat x)$ and fix a minimum point $x_n$ of this function. 
By passing to a subsequence, we may assume that $\{x_n\}$ converges to a 
point $x_\infty\in\ol B_r(\hat x)$.
Note that, since $$\,
(u-n^{-1}\phi)(x_n)\leq (u-n^{-1}\phi)(\hat x)=0,$$
we have 
\[
u(x_\infty)\leq \liminf_{n\to \infty}u(x_n)
= \liminf_{n\to\infty}(u-n^{-1}\phi)(x_n)\leq 0,
\]
which shows that $x_\infty\in\gO\stm U$,
and  that, since $u\geq 0$ in $\gO$ and $u(\hat x)=0$,  
\[
-\phi(x_n)\leq (nu-\phi)(x_n)\leq (nu-\phi)(\hat x)=-\phi(\hat x)=0,
\]
which implies 
\[
\phi(x_\infty)\geq 0 \quad \text{ and hence } \quad x_\infty=\hat x.
\]
Moreover, we may assume by passing to a subsequence that
\[
\lim_{n\to\infty}u(x_n)=\liminf_{n\to\infty}u(x_n)=0=\phi(\hat x).
\]
The viscosity property of $u$ yields together (F4) that if $n$ is large enough, then 
\[
0\geq F(x_n,u(x_n), n^{-1}D\phi(x_n),n^{-1}D^2\phi(x_n))=h(n^{-1})
F(x_n,u(x_n),D\phi(x_n),D^2\phi(x_n)),
\]
and hence, 
\[
F(x_n,u(x_n),D\phi(x_n),D^2\phi(x_n))\leq 0.
\]
Thus, in the limit as $n\to\infty$ we have 
\[
F(\hat x, \phi(\hat x),D\phi(\hat x),D^2\phi(\hat x))\leq 0. \qedhere
\]
\eproof

In the case when $u\in\LSC(\gO)$ is a supersolution of \erf{intro1} and satisfies 
$\min_\gO u=a$ for some $a\in\R$, one may obtain a result similar to the 
theorem above by applying  the theorem with $F$ and $u$ replaced, respectively, by 
$\widetilde F$ and $\tilde u$ given by
\[
\widetilde F(x,r,p,X)=F(x,r+a,p, X)   \ \ \AND \ \ \tilde u=u-a. 
\]  

\bthm \label{FU-suf0} 
 Assume \emph{(F1)--(F3)} and \emph{(F5).}
Let $U$ be an open subset of $\gO$. 
Assume that $U$ satisfies \emph{(G$_{F,\gO}$)}. Then there exists a 
supersolution $v\in \LSC(\gO)$ of \erf{intro1} such that 
\[
\min_{\gO}v=0 \ \ \ \AND \ \ \ U=\{x\in\R^N\mid v(x)>0\}. 
\]
\ethm

\begin{remark} \label{rmk1}The theorem above and the next one are considered as the existence 
results for supersolutions of \erf{intro1}, with the homogeneous Dirichlet data 
on $\gO\cap \pl U$.  
\end{remark}

\bproof  We define $v\in\LSC(\gO)$ by
\[
v(x)=\bcases
1\ \ &\IF \ x\in U,\\
0&\text{ otherwise.}
\ecases
\]
To check the supersolution property of $v$, let $\phi\in C^2(\gO)$ and $\hat x\in\gO$ and assume that $v-\phi\geq(v-\phi)(\hat x)=0$ in $\gO$.  If $\hat x\in U$, then we have $D\phi(\hat x)=0$ and $D^2\phi(\hat x)\leq 0$ by the elementary maximum principle, and moreover, by (F5), (F2) and (F3) 
\[
F(\hat x,v(\hat x),D\phi(\hat x),D^2\phi(\hat x))
\leq F(\hat x,0,0,D^2\phi(\hat x))\leq F(\hat x,0,0,0)\leq0.
\]
Otherwise, we have $\hat x\in\gO\stm U$ and 
\[
\phi(x)=-(v-\phi)(x)\leq -(v-\phi(\hat x))=\phi(\hat x)=0 \ \ \FORALL x\in\gO\stm U,
\]
and hence, by \erf{FOmU},
\[
0\geq F(\hat x, \phi(\hat x),D\phi(\hat x), D^2\phi(\hat x))
=F(\hat x, v(\hat x), D\phi(\hat x), D^2\phi(\hat x)). 
\]
This completes the proof. 
\eproof

The following theorem is a version of Theorem \ref{FU-suf0} in the 
class of continuous supersolutions.  

\bthm \label{FU-suf} Assume \emph{(F1), (F2), (F5)} and that $F=F(u,p,X)$ 
does not depend on $x$.  
Let $U$ be an open subset of $\R^N$. 
Assume that $U$ satisfies \emph{(G$_{F}$)}. Then there exists a 
supersolution $v\in \C(\R^N)$ of $F(v,Dv,D^2v)=0$ in $\R^N$ such that $v\geq 0$ in $\R^N$ and 
\[
U=\{x\in\R^N\mid v(x)>0\}. 
\]
\ethm

\bproof We set 
\[
v(x)=\dist(x, \R^N\stm U) \ \ \FOR x\in\R^N.
\]
We claim that the function $v$ is a supersolution of 
$F(v,Dv,D^2v)=0$ in $\R^N$ and, to see this, we fix  $\phi\in C^2(\R^N)$ and $\hat x\in\R^N$ 
and assume that $v-\phi$ attains a minimum at $\hat x$. 
We may assume as usual that $\phi(\hat x)=0$. 
We choose $\hat y\in\R^N\stm U$ so that $v(\hat x)=|\hat x-\hat y|$ and 
note that 
\[
|\hat x-\hat y|-\phi(\hat x)\leq |x-y|-\phi(x) \ \ \FORALL x\in\R^N,\, y\in\R^N\stm U.
\]
We plug $x=y+\hat x-\hat y$ in the above inequality, to obtain
\[
|\hat x-\hat y|-\phi(\hat x)\leq |\hat x-\hat y|-\phi(y+\hat x-\hat y) \ \ \FORALL y\in\R^N\stm U.
\] 
This says that, if we set $\psi(y):=\phi(y+\hat x-\hat y)$, then 
\[
\psi(\hat y)\geq \psi(y) \ \ \FORALL y\in\R^N\stm U.
\] 
Since $\psi(\hat y)=\phi(\hat x)=0$ and $v\geq 0$ everywhere, 
property (G$_{F,U}$) yields
\[
0\geq F(\psi(\hat y), D\psi(\hat y),D^2\psi(\hat y))\geq 
F(v(\hat x), D\psi(\hat y),D^2\psi(\hat y)),
\]
which reads 
\[
F(v(\hat x),D\phi(\hat x), D^2\phi(\hat x))\leq 0. 
\]
It is obvious that $U=\{x\in\R^N\mid v(x)>0\}$. 
\eproof

\section{Operator ${\mathcal P}^-_k$} \label{OpP_k}

In this and next sections, we deal with the operators $\cP_k^-$ and $\gl_k$, 
with $k<N$. We remark that both functions $F(X)=\cP_k^-(X)$ and $F(X)=\gl_k(X)$ satisfy (F1)--(F5). 

Let $k\in\N$ be such that $k<N$. If $F\equiv {\mathcal P}^-_k$, condition  
\erf{FOmU} reads as follows: 
\beq\tag{G$_{{\mathcal P}^-_k,\gO,U}$}\label{FU2}
\left\{\begin{minipage}{0.8\textwidth} if \ $\phi\in C^2(\gO)$,\ $\hat x\in \gO\stm U$\ and \ $\phi\leq \phi(\hat x)=0$ in $\Omega\stm U$, 
then
\[
 {\mathcal P}^-_k(D^2\phi(\hat x))\leq 0,
\]
\end{minipage}
\right. \eeq
or equivalently
\begin{equation*}
\begin{minipage}{0.85\textwidth}
if \ $\phi\in C^2(\gO)$,\ $\hat x\in \gO\stm U$\ and \ $\phi\geq \phi(\hat x)=0$ in $\gO\stm U$, 
then
\ $
{\mathcal P}^+_k(D^2\phi(\hat x))\geq 0.$
\end{minipage} 
\end{equation*}

The following proposition is an immediate consequence of the following fact: 
for any $X\in\bS^N$, 
\[
\cP_{k+1}^-(X)\leq 0\;\Rightarrow\;\cP_k^-(X)\leq 0.
\]

\begin{proposition}\label{ordP_k} Let $U$ be an open subset of $\gO$ and 
let $k\in\N$ be such that $k+1<N$. If \emph{(G$_{\cP_{k+1}^-},\gO,U$)} holds, then so does \emph{(G$_{\cP_k^-},\gO, U$)}. 
\end{proposition}

\begin{remark}\label{rmk2}\rm
Let us point out that \eqref{FU2} provides a viscosity formulation of the inequality
\begin{equation*}
\gk_{N-k}+\cdots+\gk_{N-1}\geq0
\end{equation*}
where $\gk_1\leq\cdots\leq\gk_{N-1}$ denote the principal curvature of  $\partial U$. 
If the open set $U$ has a smooth boundary,  \eqref{FU2} is indeed  equivalent, as showed in the next theorem,  to $\gk_{N-k}+\cdots+\gk_{N-1}\geq0$. 
\end{remark}


\bthm\label{thcurv} Assume that $U$ is a smooth subdomain of $\gO$, that is, the relative boundary $\gO\cap \pl U$ can be represented locally as the graph of a smooth function. Let 
\[
\gk_1\leq \gk_2\leq \cdots\leq\gk_{N-1}
\]
denote the principal curvatures of the hypersurface \ $\gO\cap\pl U$. Then we have 
\beq\label{curv}
\sum_{i=1}^k\gk_{N-i}\geq 0
\eeq
if and only if \erf{FU2} holds. 
\ethm

We use here the sign convention of the principal curvatures as follows: 
if, after a rigid transformation, 
\[
U\cap B_r =\{x=(x',x_N)\in B_r\mid x_N>f(x')\} 
\]
for some $f\in C^2(B_r^{N-1})$, with $f(0)=0$ and $Df(0)=0$, and for some $r>0$ 
such that $B_r\subset \gO$, then the eigenvalues of $D^2f(0)$ are the principal curvatures at $0\in\pl U$.

\bproof We first 
prove the sufficiency of \erf{curv} for \erf{FU2}. 
Assume that \erf{curv} holds, and let $\phi\in C^2(\gO)$ and 
$\hat x\in \gO\stm U$. Assume that $\phi(\hat x)=0$ and $\phi\leq 0$ in $\gO\stm U$. By a translation, we may assume that $\hat x=0$. 
If $0\not\in\pl U$, then $0\in\INT(\gO\stm U)$ and 
$\phi$ attains a local maximum at $0$, which yields 
\[
D^2\phi(0)\leq 0,
\]
and hence 
\[
\cP_k^-(D^2\phi(0))\leq 0.
\]

We may thus assume that $0\in\pl U$ and also that $\ol B_r\subset \gO$ for some $r>0$.  
After an orthogonal transformation,  we may assume that for some  
$f\in C^2(B_r^{N-1})$, 
\beq\label{curv1}\begin{gathered}
U\cap B_r=\{x\in B_r \mid x_N>f(x')\}, \quad  f(0_{N-1})=0,\quad Df(0_{N-1})=0_{N-1}\quad \AND \\
D^2f(0_{N-1})=\diag{\gk_1,\ldots,\gk_{N-1}}.
\end{gathered}\eeq

Note by \erf{curv1} that if $0<t<r$, then $-te_N\in \gO\stm U$ and that if 
$x\in B_r$ and $x_N=f(x')$, then $x\in\gO\stm U$. By the choice of $\phi$, 
we see that $\phi(-te_N)\leq 0$ for $0<t<r$ and that the function 
$\psi(x'):=\phi(x',f(x'))$, which is defined in an neighborhood of $0_{N-1}$,   
has a maximum at $0_{N-1}$. Accordingly, we have 
\beq\label{curv2}
\phi_{x_N}(0)=\du{D\phi(0),e_N}\geq 0 
\eeq
and, after a simple manipulation, the following matrix inequality
\beq \label{curv3}
D^2\psi(0_{N-1})=\left(D^2_{ij}\phi(0)+\phi_{x_N}(0)D^2_{ij}f(0_{N-1})\right)_{i,j\leq N-1}\leq 0.
\eeq

\def\cO{\mathcal{O}}

Let $\cO(k)$ denote the set of all collections 
$\{v_1,\ldots,v_k\}$ of orthonormal vectors in $\R^N$.  
By \erf{curv1}, \erf{curv2} and \erf{curv3} we have the matrix inequality 
\[
\left(D_{ij}^2\phi(0)\right)_{i,j\leq N-1}
\leq -\phi_{x_N}(0)\diag{\gk_1,\ldots.\gk_{N-1}},
\]
from which we deduce by \erf{curv} that
\beq\label{curv4} \bald
\cP_k^{-}(D^2\phi(0))&\,=\min_{\{v_j\}\in\cO(k)}\sum_{j=1}^k\du{D^2\phi(0) v_j,v_j}
\\&\,\leq \sum_{j=1}^k\du{D^2\phi(0) e_{N-j},e_{N-j}}
\leq -\phi_{x_N}(0)\sum_{j=1}^k \gk_{N-j}\leq 0.
\eald\eeq

From this we conclude that (G$_{\cP_k^-,\gO,U}$) holds.

Next, we assume that \erf{FU2} holds.  Let $\hat x\in\gO\cap\pl U$ and 
assume as before that $\hat x=0$. 
Select $r>0$ and $f\in C^2(B_r^{N-1})$ such that $B_r\subset \gO$ and, after 
a orthogonal transformation, 
\[
U\cap B_r=\{x\in B_r\mid x_N>f(x')\}, \quad f(0_{N-1})=0, \quad Df(0_{N-1})=0_{N-1},
\] 
and $D^2f(0_{N-1})=\diag{\gk_1,\ldots,\gk_{N-1}}$. 

For $\ga>0$ we define $\phi_\ga\in C^2(B_r)$ by 
\beq\label{curv5}
\phi_\ga (x):=\exp\big[\ga(x_N-f(x'))\big]-1.
\eeq
Note that 
\[
\phi_\ga(0)=0 \quad\AND \quad \phi_\ga(x) \leq 0 \ \ 
\IF \ x\in B_r\stm U.
\] 
We may assume by modifying $\phi_\ga$ away from the origin that $\phi_\ga\in C^2(\gO)$ and 
\[
\phi_\ga(x) \leq 0 \ \ \FORALL x\in\gO\stm U.
\]
By the assumption, we have 
\beq\label{curv6}
\cP_k^-(D^2\phi_\ga(0))\leq 0,
\eeq
while a simple computation yields
\[
D^2\phi_\ga(0)=\alpha\diag{-\gk_1,\ldots,-\gk_{N-1},\ga}. 
\]
It is now clear that for sufficiently large $\ga>0$,
\[
\cP_{k}^-(D^2\phi(0)) =-\alpha\sum_{i=1}^k\gk_{N-i},
\]
from which we conclude, together with \erf{curv6}, that 
\[
\sum_{i=1}^k\gk_{N-i}\geq 0. \qedhere
\]
\eproof

\section{Operator $\lambda_k$}\label{Opl_k}

Let $k\in\N$ be such that $k<N$. The geometric condition \eqref{FOmU} reduces in this case to the condition
\beq\tag{G$_{\lambda_k,\gO, U}$}\label{FU3}
\begin{minipage}{0.83\textwidth}
if \ $\phi\in C^2(\gO)$,\ $\hat x\in \gO\stm U$\ and \ $\phi\leq \phi(\hat x)=0$\ in\ $\gO\stm U$, 
then
\ $\lambda_k(D^2\phi(\hat x))\leq 0,$
\end{minipage}
\eeq
or equivalently to the condition
\begin{equation*}
\begin{minipage}{0.87\textwidth}
if \ $\phi\in C^2(\gO)$,\ $\hat x\in \gO\stm U$ \ and \ $\phi\geq \phi(\hat x)=0$\ in\ $\gO\stm U$, 
then
\ $\lambda_{N-k+1}(D^2\phi(\hat x))\geq 0.$
\end{minipage} 
\end{equation*}

The following proposition is a direct consequence of 
the inequality $\gl_k(X)\leq \gl_{k+1}(X)$ for all $X\in\bS^N$.

\begin{proposition}\label{ordgl_k} Let $U$ be an open subset of $\gO$ and 
let $k\in\N$ be such that $k+1<N$. If \emph{(G$_{\gl_{k+1}},\gO,U$)} holds, then so does \emph{(G$_{\gl_k},\gO, U$)}. 
\end{proposition}

\begin{remark}\label{rmk4}\rm
 \eqref{FU3} is a  viscosity definition of the inequality
$$
\gk_{N-k}\geq 0
$$
when $U$ does not necessarily have smooth boundary. Such definition is justified by the following theorem, whose proof is omitted since it can be carried out along the same line of Theorem \ref{thcurv}.   
\end{remark}


\bthm\label{thmcurv}
 Assume that $U$ is a smooth subdomain of $\gO$, that is, the relative boundary 
$\gO\cap \pl U$ can be represented locally as the graph of a smooth function. Let 
\[
\gk_1\leq \gk_2\leq \cdots\leq\gk_{N-1}
\]
denote the principal curvatures of $\pl U$. Then we have 
\beq
\gk_{N-k}\geq 0
\eeq
if and only if \erf{FU3} holds. 
\ethm

\smallskip

Now we introduce another  geometric condition  which will be proved to be equivalent to (\ref{FU3}).

\smallskip

Let $\cB_{k,m}$ be the collection of all products $\ol B_a^k \tim \ol B_b^m$, where 
$a>0$ and $b>0$. 
Let $\cC_{k,m}$ be the collection of all images of  $B\in\cB_{k,m}$ 
by rigid transformations. 
That is, 
\[
\cC_{k,m}=\{z+OB\mid z\in\R^{k+m},\ B\in\cB_{k,m},\ O\in O(k+m)\},
\] 
where $O(j)$ denotes 
the space of all $j\tim j$ orthogonal matrices. When $C=z+O\big(\ol B_a^k\tim \ol B_b^m\big)$ for some 
$z\in\R^{k+m}$, $a>0$, $b>0$ and $O\in O(k+m)$, we let $\pl'C$ denote the ``lateral  boundary'' of $C$, that is, $\pl' C=z+O\big(\pl B_a^k\tim \ol B_b^m\big)$. 

Let $U$ be an open subset of $\gO$. Consider the following condition 
concerning $U$: 
\[\begin{minipage}{0.8\textwidth} for any $\,C\in \cC_{k,N-k}$, if 
$\, C\subset\gO\,$ and  $\,\pl'C \cup \INT C\subset U$, then 
$C\subset U$.
\end{minipage}\tag{C$_{k,\gO,U}$}\label{CkOmU}
\]
When this condition holds, we say also that $U$ (resp., ($(\gO,U)$) 
satisfies (C$_{k,\gO}$) (resp., (C$_{k}$)). 
When $\gO=\R^N$, the conditions \erf{CkOmU} and (C$_{k,\gO}$) 
are referred to as (C$_{k,U}$) and (C$_k$), respectively.

\medskip

The main result of this section is the following theorem which, roughly speaking, can be summarized in the equivalence (see Remark \ref{rmk4} for the definition of $\gk_k\geq0$)
\begin{equation}\label{rough}
\text{\lq\lq\ \erf{CkOmU} $\Longleftrightarrow$ $\gk_k\geq0$ \rq\rq.}
\end{equation}

More precisely we have: 

\bthm\label{thequiv} 
Let $U$ be an open subset of $\gO$ and $k\in\N$ be such that $k<N$. Then $U$ satisfies \erf{CkOmU} if and only if $U$ satisfies 
\emph{(G$_{\gl_{N-k},\gO,U}$)}.
\ethm

The following two theorems are obvious consequences of 
Theorems \ref{FU-nec}, \ref{FU-suf} and \ref{thequiv}.

\bthm \label{nec} Let $k\in \N$ be such that $k<N$ and $\gO$ be an open subset of $\R^N$. 
Let $u\in\LSC(\gO)$ be a supersolution of 
$\,\gl_{N-k}(D^2u)=0 \ \IN \gO.\,$
Assume that $\min_{\gO}u=0$. Set 
$U=\{x\in\gO\mid u(x)>0\}$.  Then 
$U$ satisfies \emph{(C$_{k,\gO}$)}. 
\ethm

\bthm \label{suf} Let $U$ be an open subset of $\R^N$ and $k\in\N$ be such that $k<N$. Assume that $U$ satisfies \emph{(C$_{k}$)}. Then there exists a 
supersolution $v\in \C(\R^N)$ of $\gl_{N-k}(D^2v)=0$ in $\R^N$ such that $v\geq 0$ in $\R^N$ and $\,
U=\{x\in\R^N\mid v(x)>0\}.$
\ethm

An obvious result of Theorem \ref{thequiv} and Proposition \ref{ordgl_k} is the next corollary. 

\begin{corollary} \label{ordCk} Let $U$ be an open subset of $\gO$. Let $k\in\N$ satisfy $k+1<N$. If 
\emph{(C$_{k,\gO,U}$)} holds, then so does \emph{(C$_{k+1,\gO,U}$)}.  
\end{corollary}

\begin{remark}\label{rkm5} 
The existence and uniqueness results for the Dirichlet problem for the operator 
$\gl_k$ have been obtained in \cite{HL1,HL2, BR}, which should be compared 
with Theorem \ref{suf}.  The results in \cite{HL1,HL2} concern with smooth domain $U$ in $\R^N$ that satisfies $\gk_{N-k}>0$ and $\gk_k>0$ on $\pl U$.    
In \cite{BR}, the authors introduce 
a generalization of the condition $\gk_j> 0$ which makes sense for non-smooth 
domains and their assumption on the domain $U$ is stated as $\gk_{N-k}>0$ and $\gk_k>0$ on $\pl U$ in the generalized sense.  Our condition (C$_k$) for $U$ 
is similar to and weaker than the inequality $\gk_k>0$ in the generalized sense of \cite{BR}. 
\end{remark}

\subsection{Proof of the equivalence of (G$_{\gl_{N-k},\gO,U}$) and \erf{CkOmU}}
Before presenting the proof of Theorem \ref{thequiv}, 
we observe that the property \erf{FOmU}
is local in the sense as stated in the next theorem.

It is convenient to apply the conditions \erf{FOmU} and \erf{CkOmU} 
for any open subset $U$ of $\R^N$: we refer to (G$_{F,\gO,\gO\cap U}$) 
and (C$_{k,\gO,\gO\cap U}$) as \erf{FOmU} and \erf{CkOmU}, respectively.

\def\gL{\Lambda}

\begin{proposition} \label{localFOmU} Let $U$ be an open subset of $\R^N$. \emph{(i)}\  If \erf{FOmU} holds, then 
\emph{(G$_{F,\gO_0,U}$)} holds for any  open subset  $\gO_0$ of $\gO$. \emph{(ii)}\ Let $\gO_\gl$, with $\gl\in\gL$, be open 
subsets of $\R^N$ and set $\gO=\bigcup_{\gl\in\gL}\gO_\gl$.  
If \emph{(G$_{F,\gO_\gl,U}$)} holds for every $\gl\in\gL$, then \erf{FOmU} holds. 
\end{proposition}
\bproof  
We first prove assertion (i).  Let $\phi\in C^2(\gO_0)$ and $\hat x\in \gO_0\stm U$. Assume that 
\[
\phi(\hat x)=0\geq \phi(x) \ \ \FORALL x\in \gO_0\stm U. 
\]
Choose $\chi\in C^2(\gO)$ so that 
\[
\supp \chi \subset \gO_0, \ \  \chi=1 \ \text{ in a neighborhood of $\hat x$,}  
\ \ \AND \ \ \chi\geq 0 \ \ \IN \gO.
\]
Since $\supp \chi\subset\gO_0$, the product $\phi\chi$ makes sense as a function 
in $C^2(\gO)$. Note that, since $\chi\geq 0$ everywhere, 
if $x\in\gO_0\stm U$, then 
$\, (\phi\chi)(x)\leq 0, $ 
and that if $\,x\in\gO\stm \gO_0$, then $\chi(x) =0\,$ and 
$\,(\phi\chi)(x)=0\,$. Also, note that $(\phi\chi)(\hat x)=\phi(\hat x)=0$, 
$D(\phi\chi)(\hat x)=D\phi(\hat x)$ and $D^2(\phi\chi)(\hat x)=D^2\phi(\hat x)$.
Thus, we have
\[
(\phi\chi)(x)\leq 0=(\phi\chi)(\hat x) \ \ \FORALL x\in\gO\stm U,  
\]
and hence, by (G$_{F,\gO,U}$), 
\[
F(\hat x, \phi(\hat x),D\phi(\hat x), D^2\phi(\hat x))=F(\hat x, \phi\chi(\hat x),D(\phi\chi)(\hat x), D^2(\phi\chi)(\hat x))
\leq 0.
\]

We now show that assertion (ii) holds.  Let $\phi\in C^2(\gO)$ and $\hat x\in \gO\stm U$. Assume that 
\[
\phi(\hat x)=0\geq \phi(x) \ \ \FORALL x\in \gO\stm U. 
\]
Choose $\gl\in\gL$ so that $\hat x\in\gO_\gl$. Set $\psi=\phi|_{\gO_\gl}$ 
and observe that $\psi\in C^2(\gO_\gl)$ and that
\[
\psi(\hat x)=0\geq \psi(x) \FORALL x\in \gO_\gl\cap(\gO\stm U)=\gO_\gl\stm U.
\]
Since (G$_{F,\gO_\gl,U}$) is valid, we have 
\[
0\geq F(\hat x,\psi(\hat x),D\psi(\hat x),D^2\psi(\hat x))
=F(\hat x,\phi(\hat x),D\phi(\hat x),D^2\phi(\hat x)).
\]
Thus, we conclude that \erf{FOmU} holds. 
\eproof

\color{black}


\bproof[Proof of Theorem \ref{thequiv}]
We first show that (G$_{\gl_{N-k}, \gO, U}$) implies (C$_{k,\gO,U}$).  
We argue by contradiction and suppose that 
(C$_{k,\gO,U}$) does not hold.  Thus there is $C\in \cC_{k,N-k}$ such that 
\beq\label{Sep8}
C\subset\gO, \ \ \pl'C \cup\INT C\subset U \ \ \AND \ \ C\not\subset U. 
\eeq
By a rigid transformation, we may assume that 
\[
C=\ol B_a^k\tim \ol B_b^{N-k}
\]
for some $a>0,\, b>0$.  It follows from \erf{Sep8} that 
\beq\label{Sep9}
\pl B_a^k \tim \ol B_b^{N-k} \subset U, \ \ B_a^k\tim B_b^{N-k}\subset U, \ \ C\stm U\not=\emptyset \ \ 
\AND  \ \ \ol C\subset \ol U. 
\eeq
Set 
\[
T:=C\stm U=\pl C\stm U,
\]
 and note by \erf{Sep9} that 
\beq\label{Sep10}
\dist(T, \pl B_a^k\tim \ol B_b^{N-k})>0. 
\eeq

We fix a function $f\in C^2([0,\,\infty))$ so that 
\[
 f(1)=1, \ \ f'(r) <0 \ \ \AND \ \ f''(r)>0 \ \ \FORALL r\geq 0.
\]
(For example, the function $f(r)=1+\log\frac{2}{1+r}$ or $f(r)=\fr 2{1+r}$ has these properties.)
Define the function $\phi\in \C(\R^N)$ by 
\[
\phi(x)=\phi(x',x''):=\min\Big\{f\left(\fr{|x'|}{a}\right), f\left(\fr{|x''|}{b}\right)\Big\}. 
\]
Note that 
\[
\phi(x)
\bcases 
>1 & \IF \ x\in \INT C, \\ 
=1 \ & \IF \ x\in\pl C, \\
<1 & \IF \ x \in\R^N\stm C. 
\ecases
\]
Since $\INT C\subset U$, we have 
\beq\label{Sep11}
\phi(x)\leq 1 \ \ \FORALL x\in \R^N\stm U, 
\eeq
and for any $x\in\R^N\stm U$,
\beq\label{Sep12}
\phi(x)=1 \ \ \text{ if and only if } \ \  x\in T.
\eeq
Note also that $\phi$ is twice continuously differentiable on the open set
\[
E:=\{x=(x',x'')\in\R^N\stm \{0\}\mid b|x'|\not=a|x''|\}. 
\]

For any $\ep>0$ we set 
\[
\phi_\ep(x)=\phi_\ep(x',x'')=\phi(x',x'')+\ep|x''|^2. 
\]
We fix a compact neighborhood $K$ of $C$.  
For each $\ep>0$, we choose a maximum point $x_\ep$
of the function $\phi_\ep$ on the set $K\stm U$.  
We may choose a sequence $\{\ep_j\}\subset (0,\,\infty)$ converging 
to zero such that for some $x_0\in K\stm U$,
\[
x_0=\lim_{j\to \infty} x_{\ep_j}.
\]  
By  \erf{Sep9}, \erf{Sep11} and \erf{Sep12}, we have 
\[
1\leq \phi_{\ep_j}(x_{\ep_j})\leq 1+\ep_j \max_{x\in K}|x''|^2,
\]
and hence
\[
\phi(x_0)=1,
\]
which shows due to \erf{Sep12} that $x_0\in T$. Since $x_0\in\pl C$, we deduce by 
\erf{Sep10} that 
\[
x_0''\in \pl B_b^{N-k} \ \ \AND \ \ x_0'\in B_a^k,
\]  
which states that
\[
|x_0''|=b \ \ \AND  \ \ |x_0'|<a. 
\]
This shows that $x_0\in E$. 

Since $x_0\in T\subset C\subset \gO$, we may choose $r>0$ so that 
\[
B_r(x_0) \subset \gO \ \ \AND \ \ \fr{|x'|}{a}<\fr{|x''|}{b} \ \ \FORALL x\in 
B_r(x_0). 
\]
The locality of (G$_{\gl_{N-k},\gO,U})$ in regard to $\gO$, implies that the condition 
(G$_{\gl_{N-k},B_r(x_0),U}$) is valid. 
We fix  $j\in\N$ large enough so that 
$x_{\ep_j}\in B_r(x_0)$ and we write $y=x_{\ep_j}$ and $\psi=\phi_{\ep_j}$ 
for simplicity. 
Note that $\psi\in C^2(B_r(x_0))$ and  
\[
\psi(x)\leq \psi(y) \ \ \FORALL x\in B_r(x_0)\stm U. 
\]
Thanks to (G$_{\gl_{N-k},B_r(x_0),U}$), with $\phi$ replaced by the function 
$x\mapsto \psi(x)-\psi(y)$, we obtain
\[
\gl_{N-k}(D^2\psi(y))\leq 0. 
\]
 Since $|x'|/a<|x''|/b$ for $x\in B_r(x_0)$, we see that $\psi(x)=f(|x''|/b)+\ep_j |x'|^2$ 
for all  $x\in B_r(x_0)$. By a direct computation, we find that 
$D^2\psi(y)$ is similar to 
\[
\operatorname{diag}(\underbrace{2\ep_j,\ldots,2\ep_j}_{\text{$k$ times}}, f''(\gth)/b^2, f'(\gth)/(b|y''|),\ldots, f'(\gth)/(b|y''|)),
\]
where $\gth:=|y''|/b$. Hence, $D^2\psi(y)$ has $N-k-1$ negative eigenvalues and $k+1$ positive eigenvalues, which yields the contradiction $\gl_{N-k}(D^2\psi(y))> 0$. 

Now, we prove that  (C$_{k,\gO,U}$) implies  (G$_{\gl_{N-k}, \gO, U}$).  
We argue by contradiction. We thus suppose that (G$_{\gl_{N-k},\gO,U}$) does not hold, that is, 
there exist $\hat x\in \gO\stm U$ and $\phi\in C^2(\gO)$ such that 
\beq\label{Sep13}
\phi(x)\leq 0=\phi(\hat x) \ \ \FORALL x\in \gO\stm U,
\eeq
and 
\beq\label{Sep14}
\gl_{N-k}(D^2\phi(\hat x))>0. 
\eeq
By a translation, we may assume that $\hat x=0$. 
By \erf{Sep14}, we see that if we write the eigenvalues of 
$D^2\phi(0)$ in the nonincreasing order as 
\[
\mu_1\geq \mu_2\geq \cdots\geq \mu_N,
\] 
then we have 
\[
\mu_1\geq \cdots\geq \mu_{k+1}=\gl_{N-k}(D^2\phi(0))>0. 
\]
Since the condition \eqref{CkOmU} is invariant under orthogonal transformations (i.e. for 
any  orthogonal matrix $O$, then  ($\text{C}_{k,O\Omega,OU}$) holds if and only if  \eqref{CkOmU} is true), eventually replacing $\phi(x)$ by $\tilde\phi(x)=\phi(Ox)$ for some $O\in O(N)$ and using the Taylor theorem, we may choose $r>0$, $\ep>0$ and $\rho>0$ 
so that $\,B_r\subset \gO\,$ and 
\[
\phi(x)\geq p\cdot x+\ep\sum_{i\leq k+1}x_i^2-\rho\sum_{i>k+1}x_i^2 \ \ \ \FORALL x\in B_r,
\]
where $p:=D\phi(0)$.  After an  orthogonal transformation in $\R^{k+1}\tim \{0_{N-k-1}\}$, we may further assume that 
\[
p=(0_k,p_{k+1},p_{k+2},\ldots, p_N) \ \ \AND \ \ p_{k+1}\geq 0. 
\]

Let $\gth\geq \rho$ be a constant to be fixed later.  
We set 
\[
\psi(x):=p\cdot x+\ep\sum_{i\leq k+1}x_i^2-\gth \sum_{i>k+1}x_i^2. 
\]
It is obvious that $\phi\geq\psi$ everywhere and hence  we see by \erf{Sep13} that   
\beq\label{Sep15}
\text{for any $x\in B_r$, if \ $\psi(x)>0$, then \ $x\in U$.}
\eeq

Let $a>0,\,b>0$ and $t\geq 0$ be constants to be fixed later and set 
\[
C_{a,b,t}=\ol B_a^k\tim B_b^{N-k}((2\gth)^{-1}p'' +t e_{k+1}''). 
\]
Observe that 
\[\bald
\psi\big(\fr{1}{2\gth}p+te_{k+1})
&\,=p\cdot \Big(\fr{1}{2\gth}p+te_{k+1}\Big)+\ep\left(\fr{p_{k+1}}{2\gth}+t\right)^2
-\gth\sum_{i>k+1}\left(\fr{p_i}{2\gth}\right)^2 
\\&\,\geq \fr{|p|^2}{2\gth}+tp_{k+1} +\ep\left(\fr{p_{k+1}}{2\gth}+t\right)^2 -\fr{|p|^2}{4\gth}
\geq \ep t^2.
\eald 
\]

We fix a constant $L>0$ so that 
\[
\sup_{B_r}|D\psi|\leq L.
\]
Observe that if \ $\fr{1}{2\gth}p+te_{k+1}\in B_r$\ and \ $x\in B_r\cap C_{a,b,t}$, then 
\[\bald
\psi(x)&\,=\ep|x'|^2 +\psi(0_k,x'')
\\&\,\geq \ep|x'|^2 +\psi\left(\fr{1}{2\gth}p+te_{k+1}\right) 
-L\left|x''-\fr{1}{2\gth}p''-te_{k+1}''\right|
\geq \ep|x'|^2 +\ep t^2 -Lb.  
\eald
\]
Hence, if we assume, in addition, either $x\in \pl' C_{a,b,t}$ and $\ep a^2>Lb$, or, 
$\ep t^2>Lb$,  then 
\beq\label{Sep16}
\psi(x)>0. 
\eeq
Note also that if \ $b> |p|/(2\gth)$, then 
\beq\label{Sep17}
0\in C_{a,b,0}. 
\eeq

At this moment, assuming that $\gth$ is given, we fix $b$ and $a$ as 
\[
b:=\fr{|p|+1}{2\gth} \ \ \AND \ \ a:=2\sqrt{\fr{Lb}{\ep}}=\sqrt{\fr{2L(|p|+1)}{\ep\gth}}. 
\]
Noting that, as $\gth\to \infty$, the constants $a$ and $b$ converge to zero,
we may fix $\gth\geq \rho$ so that 
\[
\bigcup_{0\leq t\leq a}C_{a,b,t} \subset B_r. 
\]
In particular, we have
\[
\fr{1}{2\gth}p+te_{k+1}\in B_r \ \ \FOR t\in [0,\,a].
\]
By \erf{Sep15}, \erf{Sep16} and \erf{Sep17}, we see easily that
\beq\label{Sep18}\left\{\bald\ 
&C_{a,b,a}\subset U,
\\&\pl'C_{a,b,t}\subset U \ \ \FORALL t\in[0,a],
\\& 0\in C_{a,b,0}\backslash U.
\eald
\right.\eeq
In particular, we have 
\beq\label{Sep19}
C_{a,b,a}\stm U=\emptyset \ \ \AND \ \ 
C_{a,b,0}\stm U\not=\emptyset.
\eeq

We set 
\[
S=\{t\in [0,\,a]\mid C_{a,b,t}\stm U\not=\emptyset\} \ \ 
\AND \ \ 
\tau=\sup S. 
\]
It is easily seen that the set-valued mapping: $[0,\,a]\ni t\mapsto C_{a,b,t}\subset\R^N$ 
is upper semicontinuous in the sense that if $\{(t_n,x_n)\}_{n\in\N}
\subset[0,\,a]\tim\R^N$, 
$(t_0,x_0)\in[0,\,a]\tim\R^N$, $x_n\in C_{a,b,t_n}$ for all $n\in\N$, 
and, as $n\to \infty$, $(t_n,x_n) \to (t_0,x_0)$, then we have \ 
$x_0\in C_{a,b,t_0}$. Since $U$ is open, the set-valued mapping: 
$[0,\,a]\ni t\mapsto C_{a,b,t}\stm U\subset \R^N$ is also upper 
semicontinuous. This observation assures that $\tau$ is a maximum of $S$, 
which shows together with  \erf{Sep19} that $0\leq \tau<a$ and 
\beq\label{Sep20}
C_{a,b,\tau}\stm U\not=\emptyset,
\eeq
It follows from \erf{Sep18} that 
\beq\label{Sep21}
\pl'C_{a,b,\tau}\subset U.
\eeq
We show that 
\beq\label{Sep22}
\INT C_{a,b,\tau}\subset U.
\eeq
Indeed, otherwise, there exists a point 
\[
x_0\in \INT C_{a,b,\tau}\stm U=\left(B_a^k\tim B_b^{N-k}\left(\fr{1}{2\gth}p''
+\tau e_{k+1}''\right)\right)\stm U,
\]
which readily implies that for some $\tau<t<a$, 
\[
x_0\in \left(B_a^k\tim B_b^{N-k}\left(\fr{1}{2\gth}p''
+t e_{k+1}''\right)\right)\stm U=\INT C_{a,b,t}\stm U.
\] 
This contradicts the definition of $\tau$, which assures that \erf{Sep22} is valid. 
The relations \erf{Sep20}, \erf{Sep21} and \erf{Sep22} together 
contradict (C$_{k,\gO,U}$). 
\eproof

Combining Proposition \ref{localFOmU} and Theorem \ref{thequiv}, we see immediately the following. 

\begin{corollary} 
\label{localCkOmU} Let $U$ be an open subset of $\R^N$. \emph{(i)}\  If \erf{CkOmU} holds, then 
\emph{(C$_{k,\gO_0,U}$)} holds for any  open subset  $\gO_0$ of $\gO$. \emph{(ii)}\ Let $\gO_\gl$, with $\gl\in\gL$, be open 
subsets of $\R^N$ and set $\gO=\bigcup_{\gl\in\gL}\gO_\gl$.  
If \emph{(C$_{k,\gO_\gl,U}$)} holds for every $\gl\in\gL$, then \erf{CkOmU} holds. 
\end{corollary}

\subsection{Convexity in condition  (C$_{1,U}$)}

Concerning condition \erf{CkOmU}, we show the following proposition. 

\bthm \label{conv-in-C1} Let $U$ be an open subset of $\R^N$. Then
\emph{(C$_{1,U}$)} holds if and only if any connected component of $U$
is convex.  
\ethm

An immediate consequence of Theorems \ref{nec} and \ref{conv-in-C1} is 
the following theorem. 

\bthm\label{convexity_gl_N-1}
Let $u\in\LSC(\R^N)$ be a supersolution of $\gl_{N-1}(D^2u)=0$ 
in $\R^N\,$ and satisfy $\,\min_{\R^N}u=0$. Then any connected component of 
the positivity set $\{x\in\R^N\mid u(x)>0\}$ is convex. 
\ethm

\bproof[Proof of Theorem \ref{conv-in-C1}]  We assume first that every connected component of $U$ is convex and prove that (C$_{1,U}$) holds. To see this, we let $C\in\cC_{1,N-1}$ 
and observe that the convex hull, $\operatorname{conv} \pl'C$, of $\pl' C$ 
equals $C$. Assume that $\pl'C\cup\INT C\subset U$. Since 
$\pl'C\cup \INT C$ is connected, there is a connected component $U_0$ of $U$ 
such that $\pl'C\subset \pl'C\cup\INT C\subset U_0$.  By the convexity of $U_0$, 
we deduce that $\operatorname{conv} \pl'C\subset U_0$, which implies that 
$C\subset U_0\subset U$, which shows that (C$_{1,U}$) holds. 

We next assume that (C$_{1,U}$) holds, and prove that any connected component of $U$ is convex. 
Fix a connected component $U_0$ of $U$ and  two points $x,\,y\in U_0$. 
It is enough to prove that 
\beq \label{conv1}
[x,\,y]\subset U_0,
\eeq
where $[x,\,y]$ denotes 
the line segment 
$\{(1-t)x+ty\mid 0\leq t\leq 1\}$ connecting $x$ and $y$.  
For the proof we may assume that $x\not=y$. 

Since $U_0$ is a connected open subset of $\R^N$, 
the two points $x$ and $y$ are connected by a polygonal line in $U_0$. That is, 
there exists a finite 
collection $\{z_0,\ldots, z_m\}\subset U_0$ such that 
$x=z_0$, $y=z_m$ and $[z_{i-1},\,z_i]\subset U_0$ for all $i\in\{1,\ldots,m\}$. 
Here, we may assume that $z_{i-1}\not=z_i$ for all $i\in\{1,\ldots,m\}$. 
If $m=1$, then we have nothing to prove. 

Next we consider the case where $m\geq 2$.  To prove \erf{conv1}, we first show that   
\beq\label{conv11}
[x,\,z_2]=[z_0,\,z_2] \subset U_0. 
\eeq
If $[z_0,\,z_1]\cup [z_1,\,z_2]$ is 
a line segment, then $[z_0,\,z_2]\subset [z_0,\,z_1]\cup [z_1,\,z_2]
\subset U_0$, 
which shows that \erf{conv11} holds.  

Otherwise, the three points $z_0,\,z_1,\,z_2$
are in general position, in other words, the vectors $z_0-z_1$ and $z_2-z_1$ are 
linearly independent.  
Since $[z_0,\,z_1]\cup [z_1,\,z_2]$ is a compact subset of $U_0$, there is a 
positive number $\ep>0$ so that 
\beq\label{conv2}
B_{\ep}+[z_0,\,z_1]\cup[z_1,z_2] \subset U_0.  
\eeq
To show \erf{conv11}, we suppose to the contrary that $[z_0,\,z_2]\not\subset U_0$. 
Consider the family of triangles
\[  
\triangle_r =\{z_1 +s(z_0-z_1)+t(z_2-z_1)\mid s\geq 0,\,t\geq 0,\, s+t\leq r
\}, \ \ \text{ with } r\in(0,\,1]. 
\]
Note that $\triangle_r$ is the triangle with vertices $z_1$, $rz_0+(1-r)z_1$ and 
$rz_2+(1-r)z_1$.  We set $z_{0,r}=rz_0+(1-r)z_1$ and $z_{2,r}=rz_2+(1-r)z_0$. 
Observe by  \erf{conv2}
that if $r>0$ is sufficiently small, then $\triangle_r\subset U_0$. 
On the other hand, by the supposition that $[z_0,\,z_2]\not\subset U_0$, 
we have $\triangle_1\not\subset U_0$.   
Thus, setting 
\[
\rho:=\sup\{r\in[0,\,1]\mid \triangle_r\subset U_0\},
\] 
we have $0<\rho\leq1$. Moreover, noting that $\triangle_r\subset \triangle_t$ 
for $0<r<t\leq 1$, and that $U_0$ is open, 
we infer that 
$\triangle_\rho\not\subset U_0$ 
and $\triangle_r\subset U_0$ for all $r\in(0,\,\rho)$.  
Note also that 
\[
\bigcup_{0<t<\rho}\triangle_t \cup [z_{0,\rho},\,z_{2,\rho}]= \triangle_\rho.
\]
Thus we conclude that 
\beq\label{conv3}
[z_{0,r},\,z_{2,r}]\subset U_0 \ \ \FORALL r\in(0,\rho) 
\ \ \ \AND \ \ \ [z_{0,\rho},\,z_{2,\rho}]\not\subset U_0.
\eeq

It is clear that 
\[
\lim_{r\to\rho-} \dH([z_{0,r},\,z_{2,r}],\,[z_{0,\rho},\,z_{2,\rho}])=0,
\]
where $\dH(A,B):= 
\max\{ \sup_{x\in A}\dist(x,B),\,\sup_{y\in B}\dist(y,A)\}$ for $A, B\subset \R^N$. 
Accordingly, we may choose $r\in(0,\,\rho)$ so that 
\[
\dH([z_{0,r},\,z_{2,r}],\,[z_{0,\rho},\,z_{2,\rho}])<\ep,
\]
and, moreover, we may choose $\gd\in(0,\ep)$ so that 
\[
[z_{0,\rho},\,z_{2,\rho}] \subset B_\gd+[z_{0,r},\,z_{2,r}], 
\]
which implies together with \erf{conv3} that 
\beq\label{conv4}
(B_\gd+[z_{0,r},\,z_{2,r}])\ \stm\ U_0\not=\emptyset.
\eeq

By the compactness of $[z_{0,r},\, z_{2,r}]$ and \erf{conv3}, we deduce that 
if $t>0$ is sufficiently small, then 
\[
B_t+[z_{0,r},\,z_{2,r}] \subset U_0. 
\]
This and \erf{conv4} assure that there is $\gamma\in(0,\gd]$ such that
\beq\label{conv5}
(\ol B_\gamma+[z_{0,r},\,z_{2,r}])\stm U_0\not=\emptyset \ \ \AND \ \ 
B_\gamma+[z_{0,r},\,z_{2,r}]\subset U_0. 
\eeq
Indeed, to see this, one just needs to set
\[
\gamma=\sup\{t\in(0,\,\gd]\mid B_t+[z_{0,r},\,z_{2,r}]\subset U_0 \}.
\]

By a rigid transformation, we may assume that 
\[
z_{0,r}=0 \quad \AND \quad z_{2,r}=a e_1 \ \ \text{ for some }\ a>0.
\]

We define $C\in \cC_{1,N-1}$ by
\[
C=[0,\,a]\tim\ol B_\gamma^{N-1}.
\]
Observe by \erf{conv5} that 
\[
\INT C\subset [0,\,ae_1]+B_\gamma \subset U_0.
\]
Since
\[
B_\ep=z_{0,r}+B_\ep \subset U_0\quad \AND\quad ae_1+B_\ep=z_{2,r}+B_\ep\subset U_0\quad\text{ by \erf{conv2}},
\]
we infer from \erf{conv5} that
\[
C\stm U_0\not=\emptyset.
\]
Finally note that 
\[
\pl' C=\{0,\,ae_1\}\tim B_\gamma^{N-1}\subset B_\ep\cup (ae_1+B_\ep) \subset U_0.
\]
These together contradict (C$_{1,U}$), which ensures that \erf{conv11} holds. 

If $m=2$, then we are done. If $m>2$, then, thanks to \erf{conv11}, the collection 
of the line segments $[z_0,\,z_2], [z_2,\,z_3],\ldots,[z_{m-1},\,z_m]$ gives a polygonal line connecting $x$ and $y$ in $U_0$. Repeating this argument, we finally find that \erf{conv1} holds. 
The proof is complete. 
\eproof

\subsection{Counterexamples} The convexity  property expressed by Theorem \ref{convexity_gl_N-1} is no longer valid for supersolutions of $\lambda_k(D^2u)=0$ with $k<N-1$. Here below are counterexamples in which neither the positivity sets nor their complements are convex sets. 

\medskip

\noindent
\textbf{1. Solid torus in $\R^3$.}  
Let $V$ be a two-dimensional 
subspace of $\R^3$, $D$ be an open disk in $V$ and $l$ be a line in $V$. Assume that $l$ does not intersect $D$.  Let $U$ be an open subset 
of $\R^3$ obtained by revolving the disk $D$ about the axis $l$.  
It is easily seen that $\pl U$ has at least a positive principal curvature 
everywhere.  Theorem \ref{thmcurv} thus ensures that (G$_{\gl_1,U}$) holds, while 
$U$ and $\R^3\stm U$ are both not convex. 

\medskip

\noindent 
\textbf{2. Union of balls.} Let $B_1, B_2$ be open balls in $\R^3$ such that 
\[
B_1\stm B_2\not=\emptyset,\quad B_2\stm B_1\not=\emptyset \quad\AND\quad 
B_1\cap B_2\not=\emptyset.
\]
Set 
\[
U:=B_1\cup B_2\quad\AND \quad C:=\pl B_1\cap \pl B_2.
\]
Note that $U$ is not convex, $C$ is a circle and
\[
\pl U=\pl B_1\stm \ol B_2 \ \cup\ \pl B_2\stm \ol B_1\ \cup\ C.
\]
Observe also that for any $x\in\pl U$, there exists a unit vector $e\in\R^3$ such that 
\[
x+te \in \R^3\stm U \ \ \FORALL t\in\R.
\]
This is indeed valid for all $x\in\R^3\stm U$. 


To check that (G$_{\gl_1,U}$) is valid, let $\phi\in C^2(\R^3)$ and $\hat x\in\R^3\stm U$, and assume that $\phi\leq\phi(\hat x)=0$ in $\R^3$. We choose a unit vector 
$e\in\R^3$ so that $\hat x+te\in\R^3\stm U$ for all $t\in\R^3$. The elementary maximum principle assures that \ $\du{D^2\phi(\hat x)e,e}\leq 0,$ which implies that 
$\gl_1(D^2\phi(\hat x))\leq 0$. Thus, (G$_{\gl_1,U}$) is satisfied. 

In both cases of $U$ given above, according to Theorems \ref{FU-nec} and \ref{thequiv}, if $u$ is a supersolution of $\gl_1(D^2u)=0$ in 
$\R^3$, with the properties that $\min_{\R^3} u=0$ and $U=\{x\in\R^3\mid  u(x)>0\}$, then (C$_{2,U}$) is valid.

\section{The strong maximum principle}\label{smp}

In this section, we are concerned with the strong maximum principle 
for degenerate elliptic equation \erf{intro4}.  
Our original interest  in this matter is for the operators $\cP_k^-$ and $\gl_k$. 

In this section, let $\Omega$ be an open subset of $\mathbb R^N$, as before, 
and let $f\in \C(\Omega)$ satisfy $f\leq0$ in $\gO$. 
The strong maximum principle does not hold in general, as the examples below show. 
Extra assumptions on $f$ are needed for the validity of \eqref{SMP}. It is somehow related to the geometry of the set $$\Gamma=\left\{x\in\Omega:\,f(x)=0\right\}.$$ 

A simple and useful observation is 
that if $u$ satisfies the assumptions of \erf{SMP} and 
attains a minimum value zero at $x_0\in\gO$, then we have $0=F(x_0,0,0,0)\leq f(x_0)$, and hence, $x_0\in \Gamma$, which implies that
\beq\label{smp2}
\{x\in\gO\mid u(x)=0\}\subset \Gamma.
\eeq

Note that in the next examples we consider $F=\lambda_k$, since it is clear that any supersolution of $\lambda_k(D^2u)=0$ is automatically a supersolution of ${\mathcal P}^-_k(D^2u)=0$. The examples tell us not only that the degeneration of ellipticity 
of the operator $\gl_k$ violates \erf{SMP}, but also that the $N$--dimensional 
Lebesgue measure of $\Gamma$, a geometric quantity, 
does not give a criterion for the validity of 
\erf{SMP}. 

\begin{example}\rm
Let $B_R\Subset\Omega$ and let $\alpha>2$. The function $$u(x)=\begin{cases}
{(R^2-|x|^2)}^\alpha & \text{if $|x|<R$}\\
0 & \text{otherwise}
\end{cases}$$
is a $C^2$-solution of the equation 
$$
\lambda_k(D^2u)=f(x)\qquad\text{in\; $\Omega$}
$$ 
for  $k<N$, where 
$$f(x)=\begin{cases}
{-2\alpha(R^2-|x|^2)}^{\alpha-1} & \text{if $|x|<R$}\\
0 & \text{otherwise.}
\end{cases}$$
This contradicts \eqref{SMP}. Note that $\Gamma=\Omega\backslash B_R$, then $|\Gamma|$, the $N$--dimensional Lebesgue measure of  $\Gamma$, is positive. 

This is also a consequence of the general results in the previous sections.  
Since $B_R$ is a convex set, by Theorem \ref{conv-in-C1}, (C$_{1,B_R}$) holds
and, consequently, thanks to Theorems \ref{FU-suf} and \ref{thequiv}, there exists a supersolution $v\in \C(\R^N)$ of 
$\gl_{N-1}(D^2v)=0$ in $\R^N$ such that $B_R=\{x\in\R^3\mid v(x)>0\}$.
Note here that $v$ is a supersolution of $\gl_k(D^2v)=0$ in $\gO$ for any $k<N$. 
According to \erf{smp2}, we have $\gO\stm B_R\subset \Gamma$ and  $|\Gamma|$ is positive. 
\end{example}

\smallskip
Next example shows that not even the condition $|\Gamma|=0$ is sufficient to ensure \eqref{SMP}. 

\begin{example}\rm
Let $\Omega\subset\mathbb R^N$ be an open set such that $0\in\Omega$. Let $k<N$ and set 
\[
U=\{x\in\R^N\mid \prod_{i=1}^k x_i\not=0\}. 
\]
This $U$ can be represented as the direct sum of its connected components:
\[
U=\bigcup_{\gs} U_\gs,
\]
where $\gs$ ranges over all $\gs=(\gs_1,\ldots,\gs_k)$, with $\gs_i\in\{-1,\,1\}$, 
and 
\[
U_\gs=\{x\in\R^N\mid \gs_1 x_1>0,\ldots, \gs_k x_k>0\}. 
\]

Assume $N\geq2k$ for some $k\in\mathbb N$. Consider the function 
$$u(x)=\sum_{i=1}^k|x_i|^\alpha$$ 
where $\alpha\in(0,1)$ and $x=(x_1,\ldots,x_N)\in\Omega$.
We have
$$
D^2u(x)=\text{diag}\left[\alpha(\alpha-1)|x_1|^{\alpha-2},\ldots,\alpha(\alpha-1)|x_k|^{\alpha-2},0,\ldots,0\right]\qquad\text{if}\quad\prod_{i=1}^kx_i\neq0.
$$
Then
\begin{equation}\label{smp3}
\lambda_k(D^2u)=\max\left\{\alpha(\alpha-1)|x_1|^{\alpha-2},\ldots,\alpha(\alpha-1)|x_k|^{\alpha-2}\right\}
\end{equation}
whenever $x\in\Omega$ and $x_i\neq0$ for any $i=1,\ldots,k$.

We claim that if $x_i=0$ for some $i=1,\ldots,k$ then the function $u$ satisfies the inequality
\begin{equation}\label{smp4}
\lambda_k(D^2u(x))\leq0
\end{equation}
in the viscosity sense. For this let  $\varphi\in C^2(\Omega)$ be such that 
\begin{equation}\label{smp5}
\min_{\Omega}(u-\varphi)=(u-\varphi)(x)
\end{equation}
for some $x\in\Omega$ with $x_i=0$.

We choose a $k$-dimensional subspace $V_k$ of $\mathbb R^N$ orthogonal to $e_1,\ldots,e_k$. Note that such choice is possible since $N\geq2k$. In this way for any unit vector $v\in V_k$ and any $t\in\mathbb R$ such that $x+tv\in\Omega$ we have
$$
u(x+tv)=u(x).
$$
From \eqref{smp5} we have 
$$
\varphi(x+tv)\leq\varphi(x)
$$ 
and then 
$$
\left\langle D^2\varphi(x)v,v\right\rangle\leq0.
$$
From the  variational characterization of the
eigenvalues 
\begin{equation*}
\lambda_{k}(D^2\varphi(x_0))=\min_{\dim V=k}\max_{v\in V,\,|v|=1}\left\langle D^2\varphi(x)v,v\right\rangle
\end{equation*}
we deduce that $\lambda_{k}(D^2\varphi(x))\leq0$, as claimed.

Putting together \eqref{smp3} and \eqref{smp4}, we infer that $u(x)$ is a viscosity supersolution of 
$$
\lambda_k(D^2u)=f(x)\qquad\text{in $\Omega$},
$$
where $f$ is any continuous function such that
$$
0\geq f(x)\geq\max\left\{\alpha(\alpha-1)|x_1|^{\alpha-2},\ldots,\alpha(\alpha-1)|x_k|^{\alpha-2}\right\}  \qquad\text{if\;\; $\prod_{i=1}^kx_i\neq0$}
$$
and 
$$
f(x)=0 \qquad\text{if \;\;$\prod_{i=1}^kx_i=0$}.
$$
In this case $\Gamma$ is given by the union of $k$ hyperplanes 
$\{x\in\R^N\mid x_i=0\}$, with $i=1,\ldots,k$, and, hence $|\Gamma|=0$.
\end{example}

We give  a general sufficient condition on $\Gamma$ for the validity of \eqref{SMP}. 

\bthm \label{SMP-suff} Assume \emph{(F1), (F2)} and \emph{(F4)} and that 
$\Gamma\not=\gO$. Let  $u\in \LSC(\Omega)$ be a supersolution of  
\begin{equation}\label{smp6}
F(D^2u)=f\qquad\text{in\; $\Omega$.}
\end{equation} 
Assume $u\geq 0$ in $\gO$ and that for any nonempty closed subset $\Gamma_0$ of $\Gamma$, condition \emph{(G$_{F,\gO,U}$)}, with $U:=\gO\stm \Gamma_0$, does not hold. Then $u>0$ in $\Omega$.
\ethm

\begin{remark} The conclusion of the theorems in this section is stronger than the 
standard strong maximum principle in that there is no alternative conclusion 
of supersolutions being identically zero.   
\end{remark}

\bproof  We suppose by contradiction that $\min_\gO u=0$.  
Set $U=\{x\in\gO\mid u(x)>0\}$ and $\Gamma_0=\gO\stm U$. 
Note that $\Gamma_0=\{x\in\gO\mid u(x)=0\}$ is a nonempty closed (in the relative topology) subset of 
$\gO$.   
By \erf{smp2}, we have $\Gamma_0\subset \Gamma$. 
By the assumption, condition (G$_{F,\gO,U}$) does not hold, which contradicts 
 Theorem \ref{FU-nec}.  This contradiction completes the proof.  
\eproof 

The same proof as above ensures the following corollary.

\begin{corollary}\label{SMP-suff-cor} Assume \emph{(F1), (F2)} and \emph{(F4)} and that 
$\Gamma\not=\gO$. Let  $u\in \LSC(\Omega)$ be a supersolution of  
\erf{smp6}. 
Assume $u\geq 0$ in $\gO$, that there is a closed subset $\Gamma'$ of $ \Gamma$ such that $\{x\in\gO\mid u(x)=0\} \subset \Gamma'$ and
that for any nonempty closed subset $\Gamma_0$ of $\Gamma'$, condition \emph{(G$_{F,\gO,U}$)}, with $U:=\gO\stm \Gamma_0$, does not hold. Then $u>0$ in $\Omega$. 
\end{corollary}

Henceforth, $F$ is either $\lambda_k$ or ${\mathcal P}^-_k$.  
The following theorem has a little more explicit conditions for \erf{SMP} to hold.

\begin{theorem}\label{SMP-suff2}
Let  $u\in \LSC(\overline\Omega)$ be, as in the previous theorem, 
a supersolution of  \erf{smp6} and satisfies $u\geq 0$ in $\gO$. 
If one of the following conditions holds, then $u$  cannot achieve its minimum inside $\Omega$:
\begin{itemize}
	\item[(i)] $\Gamma$ consists of isolated points.  
	\item[(ii)] There exists $\delta>0$ such that $$\Gamma\subset\Omega_\delta:=\left\{x\in\Omega:\,{\rm dist}(x,\partial\Omega)>\delta\right\}.$$
	\item[(iii)]  Let $F=\lambda_k$ and assume that for any $x\in\Gamma$ there exist positive numbers $r_1=r_1(x),\,r_2=r_2(x)$ and an integer $1\leq h=h(x)\leq k-1$ such that, up to an orthogonal transformation around $x$, 
\[\bcases
	\left(\,\overline B_{r_1}^h(x_1,\ldots,x_h)\times\partial B_{r_2}^{N-h}(x_{h+1},\ldots,x_N)\right)\cap \Gamma=\emptyset, &\\[3pt]
\ol B_{r_1}^h(x_1,\ldots,x_h)\tim \ol B_{r_2}^{N-h}(x_{h+1},\ldots,x_N)\subset \gO.&
\ecases
\]
	\item[(iv)] Let $F={\mathcal P}^-_k$ and assume that for any $x\in\Gamma$ there exist positive numbers $r_1=r_1(x),\,r_2=r_2(x)$ and an integer $1\leq h=h(x)\leq k-1$ such that 
	$\left(\frac{r_2}{r_1}\right)^2<\frac kh-1$ and
	, up to an orthogonal transformation around $x$, 
\[\bcases
	\left(\,\overline B_{r_1}^h(x_1,\ldots,x_h)\times\partial B_{r_2}^{N-h}(x_{h+1},\ldots,x_N)\right)\cap \Gamma=\emptyset,&\\[3pt]
\ol B_{r_1}^h(x_1,\ldots,x_h)\tim \ol B_{r_2}^{N-h}(x_{h+1},\ldots,x_N)\subset \gO.&
\ecases\]
\end{itemize}
\end{theorem}

The following proof is based on the comparison principle for \erf{smp6}, with 
$F=\gl_k$ or $F=\cP_k^-$. The validity of this maximum principle  
is guaranteed by observing that $\gl_k(X+t I)=\gl_k(X)+t$ and 
$\cP_k^-(X+tI)=\cP_k^-(X)+kt$ for any $X\in\bS^N$ and $t\in\R$, which implies that 
if $u$ is a supersolution of \erf{smp6}, then the  function $u(x)-\ep|x|^2$ is a 
strict supersolution of \erf{smp6}, and by recalling the general strategy 
explained in \cite[5.C]{CIL}. See also \cite[Proposition~3.3]{GV} and \cite{HL1}

\begin{proof}
For the proof, we argue by contradiction in any case and suppose 
that $u$ attains its minimum at $x_0\in\Omega$. Without loss of generality may assume $x_0=0$ and $u(0)=0$. Since $F(D^2u)\leq f$ in $\Omega$ then necessarily $0\in\Gamma$. 

Consider first the case where condition (i) holds. 
Since $0$ is an isolated of $\Gamma$,  there is $\delta>0$ such that $\overline B_\delta\subset\Omega$ and 
\begin{equation*}
\overline B_\delta\cap\Gamma=\left\{0\right\}.
\end{equation*} 
Since $f<0$ on ${\partial B_\delta}$ we infer that 
\begin{equation}\label{smp7}
\min_{\partial B_\delta}u>0.
\end{equation}
 On the other hand it is standard to prove that the operator $F$ satisfies the comparison principle, hence $\displaystyle\min_{\overline B_\delta}u=0=\min_{\partial B_\delta}u$ so contradicting \eqref{smp7}. 

\medskip

In a similar way we get a contradiction under the assumption (ii). Indeed, since $f<0$ in $\partial\Omega_{\delta}$ we have
$$
\min_{\partial\Omega_{\delta}}u>0.
$$
But this is impossible, again by the comparison principle. 

\medskip
Assume now that condition (iii) holds. Let 
\begin{equation}\label{smp8}
\varphi(x)=\alpha\left[\,\sum_{i=h+1}^Nx_i^2-\beta\sum_{i=1}^hx_i^2\,\right],
\end{equation}
where $\alpha,\beta$ are positive constant to be fixed in such a way $u-\varphi$ has a local minimum point inside $B_{r_1}^{h}\times B_{r_2}^{N-h}$.\\
Fix $\beta>\left(\frac{r_2}{r_1}\right)^2$. For any $x\in\partial B_{r_1}^h\times \overline B_{r_2}^{N-h}$ we have
\begin{equation}\label{smp9}
(u-\varphi)(x)\geq -\varphi(x)\geq\alpha\left[\,\beta r_1^2-r_2^2\,\right]>0.
\end{equation} 
Now we choose $\alpha$ positive and small enough such that 
\begin{equation}\label{smp10}
\min_{\overline B_{r_1}^h\times\partial B_{r_2}^{N-h}}u>\alpha r_2^2.
\end{equation}
Note that such choice  is possible since, by assumption iii), $\overline B_{r_1}^h\times\partial B_{r_2}^{N-h}$ is a compact subset of $\left\{f<0\right\}$ where $u>0$. By \eqref{smp10}, for any $x\in \overline B_{r_1}^h\times\partial B_{r_2}^{N-h}$ we have
\begin{equation}\label{smp11}
(u-\varphi)(x)=u(x)-\alpha\left[\,r_2^2-\beta\sum_{i=1}^hx_i^2\,\right]\geq \min_{\overline B_{r_1}^h\times\partial B_{r_2}^{N-h}}u-\alpha r_2^2>0.
\end{equation}
Since $u(0)=\varphi(0)$, from \eqref{smp9} and \eqref{smp11} we infer that there exists $\hat x\in B_{r_1}^h\times B_{r_2}^{N-h}$ such that 
$$
(u-\varphi)(\hat x)=\min_{B_{r_1}^h\times B_{r_2}^{N-h}}(u-\varphi)\,.
$$
Using the equation \eqref{smp6} we obtain a contradiction:
\begin{equation*}
\begin{split}
0\geq f(\hat x)
&\geq \lambda_k(D^2\varphi(\hat x))\\&=2\alpha\,\lambda_k({\rm diag}[\underbrace{-\beta,\ldots,-\beta}_{\text{$h$  times}},\underbrace{1,\ldots,1}_{\text{$N-h$ times}}])\\
&=2\alpha>0.
\end{split}
\end{equation*}

\smallskip

Under the assumption (iv), we still consider the test function \eqref{smp8} and by the argument above we get
\begin{equation*}
\begin{split}
0\geq f(\hat x)
&\geq {\mathcal P}^-_k(D^2\varphi(\hat x))\\&=2\alpha\,{\mathcal P}^-_k({\rm diag}[\underbrace{-\beta,\ldots,-\beta}_{\text{$h$  times}},\underbrace{1,\ldots,1}_{\text{$N-h$ times}}])\\
&=2\alpha(-h\beta+k-h).
\end{split}
\end{equation*}
This again leads to a contradiction if $\beta<\frac kh-1$.
 \end{proof}

The proof above is mainly based on the comparison principle for the operators 
$\cP_k^-$ and $\gl_k$.  The following proof of Theorem \ref{SMP-suff2}, (ii) 
is based on Theorem \ref{SMP-suff}. Notice that, in Theorem \ref{SMP-suff2}, 
claim (i) is an easy corollary of (ii).   

\bproof[Proof of Theorem \ref{SMP-suff2}, (ii)] Suppose to the contrary that 
$\min_\gO u=0$ 
We choose a constant $M>0$ so that $M<\min_{\pl\gO_\gd}u$ and 
set
\[
v(x)= \bcases
\min\{M,u(x)\} \ \ & \FOR \ x\in\gO_\gd, \\ 
M&\FOR \ x\in \R^N\stm \gO_\gd. 
\ecases
\]
Observe that $v\in\LSC(\R^N)$ and it is a supersolution of $\gl_k(D^2v)=0$ in $\R^N$. Moreover, $v$ is a supersolution of $\gl_1(D^2v)=0$ in $\R^N$. 
Note also that   
\[
\{x\in\R^N\mid v(x)=0\}=\{x\in\gO\mid u(x)=0\}\subset \Gamma. 
\]

Thanks to Corollary \ref{SMP-suff-cor}, we only need to show that 
if $\Gamma_0$ is a closed subset of $\Gamma$, then (G$_{\gl_1,U}$), with 
$U:=\R^N\stm \Gamma_0$, does not hold. Recall by Theorem \ref{thequiv} that 
 (G$_{\gl_1,U}$) is equivalent to condition (C$_{N-1,U}$).

We choose a constant $R>0$ so that $\Gamma_0\subset B_R$.  
Let $b>0$ and consider 
\[
C_b=\ol B_R^{N-1}\tim [-R,-R+b] \in\cC_{N-1,1}.
\]
It is clear that $\pl' C_b=\pl B_R^{N-1}\tim[-R,-R+b] \subset U$. 
Observe also that $C_b\cap \Gamma_0=\emptyset$ for sufficiently small 
$b>0$ and that $C_b\cap \Gamma_0\not=\emptyset$ for sufficiently large $b>0$. 
We may select $b>0$ so that $C_b\cap \Gamma_0\not=\emptyset$ 
and $\INT C_b\cap\Gamma_0=\emptyset$. 
Thus, we see that $\pl' C_b\cup\INT C_b\subset U$ and $C_b\not\subset U$, which 
means that (C$_{N-1,U}$) does not hold. The proof is complete. 
\eproof

In the same spirit as above we give another proof of Theorem \ref{SMP-suff2}, (iii).

\bproof[Proof of Theorem \ref{SMP-suff2}, (iii)]
We suppose by contradiction that $\min_{\gO} u=0$ and set $U=\{x\in\gO\mid u(x)>0\}$ and $\Gamma_0=\gO\stm U$.   
As before, we infer that  $\Gamma_0\subset\Gamma$ and that $U$ satisfies (C$_{N-k,\gO,U}$).
Moreover, according to Corollary \ref{ordCk}, condition (C$_{N-h-1,\gO,U}$) holds. 

Noting that $\Gamma_0\not=\emptyset$, we fix $\hat x\in \Gamma_0$. 
By the assumption, we may assume that there exist $a>0$ and $b>0$ 
such that  
\[
C:=\ol B_a^{N-h}(\hat x_1,\ldots,\hat x_{N-h})\tim \ol B_b^{h}(\hat x_{N-h+1},\ldots,\hat x_N) 
\subset \gO, 
\]
and
\[
\pl' C=\pl B_a(\hat x_1,\ldots,\hat x_{N-h})\tim \ol B_b(\hat x_{N-h+1},\ldots,\hat x_N)\cap \Gamma=\emptyset.
\]
By a translation, we may assume that $\hat x=0$. 

Fix a constant $\gl>0$ so that 
\[
\gl b^2>a,
\]
and consider the function 
\[
g(x):=|x'|-\gl|x''|^2 \ \ \ON \R^N,
\]
where $x'=(x_1,\ldots,x_{N-h})\in\R^{N-h}$ and $x''=(x_{N-h+1},\ldots,x_N)\in\R^{h}$.  
Set 
\[
m:=\max_{C\stm U} g,
\]
and note that
 if \ $x\in C$\ and \ $g(x)>m$, then \ $x\in U$.
Since $\hat x=0\in C\stm U$, we have $m\geq 0$. 
Let $x\in \pl C$. If $|x'|=a$, then $x\in\pl'C\subset U$. 
Otherwise, we have $|x''|=b$ and 
\[
g(x)\leq a -\gl b^2<0\leq m.
\]

We fix a maximum point $\xi\in C\stm U$ of $g$ on $\gO\stm U$. 
The observations above show that 
\beq \label{smp12} \xi\in \INT C.\eeq 
After an orthogonal transformation in $\R^{N-h}\tim \{0_{h}\}$, 
noting that $g$ is invariant under such a transformation, we may assume that 
\[
\xi'=|\xi'|\,e_{N-h}. 
\]

In what follows we consider the paraboloid 
\beq\label{smp13}
r=\gl|\eta|^2+m
\eeq
in $\R^{h+1}, $where $r\in\R$ and $\eta\in\R^{h}$ and consider balls tangent to this 
paraboloid. In view of \erf{smp12}, we have 
\[
|\xi'|<a \  \ \AND  \ \ |\xi''|<b.
\] 
Noting that $(r,\eta)=(|\xi'|,\xi'')$ satisfies \erf{smp13}, we may choose 
$\ep\in (0,\,\min\{a,\,b\})$ and $(r_0,\eta_0)\in\R^{h+1}$ so that the ball $B_\ep^{h+1}(r_0,\eta_0)$ is tangent 
to the paraboloid \erf{smp13} at $(|\xi'|,\xi'')$, it stays in the positive side with regard to the $r$-axis, and it stays ``inside $C$''. More precise requirements on the ball $B_\ep^{h+1}(r_0,\eta_0)$  are stated as 

\begin{align}
&|(|\xi'|,\xi'')-(r_0,\eta_0)|=\ep,\label{smp14}
\\& B_\ep^{h+1} (r_0,\eta_0)\subset \{(r,\eta)\in\R\tim\R^{h}\mid r> \gl |\eta|^2+m\},
\label{smp15}
\\& \ol B_{2\ep}^{h+1}(r_0,\eta_0) \subset [-a,\,a]\tim \ol B_b^{h}. 
\label{smp16}
\end{align}

We set 
\[
C_\ep=\ol B_\ep^{N-h-1}\tim \ol B_\ep^{h+1}(r_0,\eta_0) \in\cC_{N-h-1,h+1}.
\]

We claim that 
\[
C_\ep\subset C. 
\]
To see this, let $x\in C_\ep$ and write $x=(\tilde x,r,\eta)\in\R^{N-h-1}\tim \R\tim\R^{h}$. 
Note by \erf{smp16} that 
\[
r\in [r_0-\ep,r_0+\ep]\subset [-a+\ep,a-\ep],   
\]
and hence
\[
|x'| \leq |\tilde x|+|r|\leq \ep+a-\ep=a. 
\]
That is, we have \ $x'\in\ol B_a^{N-h}$, and also, by \erf{smp16} we have 
\ $\eta\in \ol B_b^{h}$.  
Thus, we see that 
\beq\label{smp17}
C_\ep \subset C\subset \gO.
\eeq

Note by \erf{smp14} that 
\beq\label{smp18}
\xi=(0_{N-h-1},|\xi'|,\xi'')\in B_\ep^{N-h-1}\tim\pl B_\ep^{h+1}(r_0,\eta_0)\subset C_\ep.
\eeq

We show next that 
\beq\label{smp19}
\INT C_\ep \subset U. 
\eeq
Let $x=(x',x'')=(\tilde x,r,x'')\in\INT C_\ep=B_\ep^{N-h-1}\tim B_\ep^{h+1}(r_0,\eta_0)$, where, as above, $x'=(\tilde x,r)\in\R^{N-h}$ 
and $x''\in\R^{h}$.  Since $(r,x'')\in B_\ep(r_0,\eta_0)$, the inclusion \erf{smp15} implies that 
\[
|x'|\geq r>\gl|x''|^2+m. 
\]
Since $x\in C$, we deduce by the choice of $m$ that $x\in U$. 
Thus, \erf{smp19} is valid. 

Next, let   $x\in\pl'C_\ep=\pl B_\ep^{N-h-1}\tim\ol B_\ep^{h+1}(r_0,\eta_0)$ 
and set $x=:(x',x'')\in \R^{N-h}\tim \R^{h}$ 
and $x'=:(\tilde x,r) \in\R^{N-h-1}\tim \R$. The inclusion 
$(r,x'')\in\ol B_\ep(r_0,\eta_0)$ and \erf{smp15} yield 
\[
r\geq \gl|x''|^2+m,
\]
which implies that
\[
|x'|=\sqrt{\ep^2+r^2}>r\geq \gl|x''|^2+m.
\]  
Since $x\in C$, this shows that $x\in U$. Hence, we have 
\[
\pl'C_\ep\subset U.
\] 
This, \erf{smp17}, \erf{smp18} and \erf{smp19} ensure that 
(C$_{N-h-1,\gO,U}$) does not hold, which is a contradiction. 
The proof is complete. 
\eproof

\end{document}